\documentclass[12pt]{amsart}
\usepackage{amssymb,amsmath}
\usepackage{epsfig}
\usepackage{graphicx}
\usepackage{colordvi}
\numberwithin{equation}{section}

\newtheorem{theorem}[equation]{Theorem}
\newtheorem{lemma}[equation]{Lemma}
\newtheorem{proposition}[equation]{Proposition}
\newtheorem{cor}[equation]{Corollary}

\theoremstyle{definition}

\newtheorem{example}[equation]{Example}

\newtheorem{remark}[equation]{Remark}

\copyrightinfo{}{}
\newcounter{FNC}[page]
\def\fauxfootnote#1{{\addtocounter{FNC}{2}$^\fnsymbol{FNC}$%
     \let\thefootnote\relax\footnotetext{$^\fnsymbol{FNC}$#1}}}

\newcommand{\QED}{\hfill$\Box$}

\newcommand{\Span}{\mbox{\rm Span}}
\newcommand{\Hom}{\mbox{\rm Hom}}
\newcommand{\Mat}{\mbox{\rm Mat}}

\newcommand{\chibar}{{\overline{\chi}}}

\newcommand{\calO}{{\mathcal O}}

\newcommand{\refl}{\mbox{\tiny $\left(\hspace{-3pt}\begin{array}{rr}-1\hspace{-5pt}&0\\
            0\hspace{-5pt}&1\end{array}\hspace{-3pt}\right)$}}

\newcommand{\iden}{\mbox{\tiny $\left(\hspace{-3pt}\begin{array}{rr}1&0\\
            0&1\end{array}\hspace{-3pt}\right)$}}

\newcommand{\C}{{\mathbb C}}
\renewcommand{\H}{{\mathbb H}}

\renewcommand{\P}{{\mathbb P}}
\newcommand{\R}{{\mathbb R}}
\begin{document}
\title{Convex Hulls of Orbits and Orientations of a Moving Protein Domain}

\author{Marco Longinetti}
\address{Dipartimento ingegneria agraria e forestale\\
         Universit\`a degli Studi di Firenze\\
         Piazzale delle Cascine, 15\\
         50144 Firenze\\
         ITALIA}
\email{longinetti@diaf.unfi.it}
\urladdr{\SMALL{wwwnt.unifi.it/diaf/nuovosito/withframes/docenti/longinetti/sitolonginetti/index.htm}}
\author{Luca Sgheri}
\address{IAC - CNR Sede di Firenze\\
         Via Madonna del Piano 10\\
         50019 Sesto Fiorentino (FI)\\
         ITALIA}
\email{luca@fi.iac.cnr.it}
\urladdr{www.fi.iac.cnr.it/iaga/luca.html}
\author{Frank Sottile}
\address{Department of Mathematics\\
         Texas A\&M University\\
         College Station\\
         TX \ 77843-3368\\
         USA}
\email{sottile@math.tamu.edu}
\urladdr{www.math.tamu.edu/\~{}sottile}
\thanks{Work of Sottile supported by NSF CAREER grant DMS-0538734 and Peter
  Gritzmann of Technische Universit\"at M\"unchen}
\begin{abstract}
 We study the facial structure and Carath\'eodory number of
 the convex hull of an orbit of the group of rotations in $\R^3$ acting on the
 space of pairs of anisotropic symmetric $3\times 3$ tensors.
 This is motivated by the problem of determining the structure of some
 proteins in aqueous solution.
\end{abstract}

\maketitle
%
%
\section*{Introduction}

The most aesthetically appealing polytopes arise as convex
hulls of orbits of finite groups acting on a vector space.
These include the platonic and archimedean solids and their higher-dimensional
generalizations, such as the regular polytopes~\cite{Coxeter}.
In contrast, the analogous objects for compact Lie groups have not
attracted much study.
We investigate convex hulls of orbits of the group
$SO(3)$ in a particular 10-dimensional representation.
Our motivation comes from an algorithm to understand the fold of some proteins.

Certain proteins, such as calmodulin~\cite{BGKLP},
consist of two rigid domains connected via a region that is flexible
in aqueous solution (i.e.~under physiological conditions), and the
problem is to determine the relative position and orientation of
these two domains. Calmodulin, as many other proteins, incorporates
metal ions into its structure. When a paramagnetic ion is
substituted, it interacts with the magnetic field of dipoles within
the protein via its magnetic susceptibility tensor $\chi$. Part of
this interaction, the residual dipolar coupling, may be inferred
from nuclear magnetic resonance data and depends solely upon the
relative orientation of the two domains. When the relative
orientation of the two domains is not constant, we infer the mean
magnetic susceptibility tensor $\chibar$ from these data.

We model this relative orientation by a probability measure $p$ on the
group $SO(3)$ of rotations of $\R^3$.
Then $\chibar$ is the average with respect to $p$ of rotations of $\chi$.
Recovering $p$ from $\chibar$ is an ill-posed inverse problem.
Nevertheless, $\chibar$ contains useful information about $p$.
Gardner, Longinetti, and Sgheri~\cite{GLS05} gave an algorithm
to determine the maximum probability of a given relative
orientation of the two domains.
Since $\chibar$ lies in the convex hull $V$ of the orbit of $\chi$ under the group
of rotations of $\R^3$, it admits a representation
$\chibar = \sum_j p_j R_j.\chi$,
%
%
where the sum is finite, $\sum_jp_j=1$ with $p_j\geq 0$, $R_j$ is a rotation in $\R^3$,
and $R_j.\chi$ is the action of $R_j$ on the tensor $\chi$.
The minimal number of
summands needed 
to represent any point $\chibar$ in $V$ is the
Carath\'eodory number of $V$.

It is often possible to substitute a different metal ion into the
protein with a different susceptibility tensor $\chi'$. Repeating
the measurements gives a second mean tensor $\chibar'$ which is the
average of rotations of $\chi'$ with respect to the measure $p$.
Combining this with $\chibar$ gives more information about $p$.
Longinetti, Luchinat, Parigi, and Sgheri~\cite{LLPS} adapt the
algorithm of~\cite{GLS05} when there are two or more metal ions and
show how this can be used to better understand the structure of
calmodulin. Their algorithm uses some knowledge of the convex hull
$V^{1,2}$ of the orbit of the pair $(\chi,\chi')$ under the group of
rotations.

We study the Carath\'eodory number and facets of $V^{1,2}$. When
$\chi$ and $\chi'$ are linearly independent, $V$ has dimension 10.
We call a subgroup of $SO(3)$ which stabilizes a line in $\R^3$ a
\Blue{{\it coaxial group}}, and a face of $V^{1,2}$ which is
stabilized by a such a subgroup a \Blue{{\it coaxial face}}. Our
main result is the following.\medskip

\noindent{\bf Theorem~\ref{T:coaxial}.}
 {\it
  Faces of $V^{1,2}$ have dimension at most $6$.
  The coaxial faces of\/ $V^{1,2}$ form a $3$-dimensional family
  whose union is a $9$-dimensional
  subset of the boundary of\/ $V^{1,2}$ if and only if $\chi$ and $\chi'$ have distinct
  eigenvectors.
  In that case, almost all coaxial faces have dimension $6$, have Carath\'eodory number
  $4$, and are facets of\/ $V^{1,2}$.  }\medskip

Our main result implies that the Carath\'eodory number of $V^{1,2}$
is at most $8$. This is an advance over \cite{GLS05}, where it has
been bounded between $4$ and $10$ inclusive.

We are unable to show that the boundary of $V^{1,2}$ is the union of
coaxial faces when $\chi$ and $\chi'$ have distinct eigenvectors,
but conjecture that this is the case. As a consequence of our main
result, we also conjecture that the Carath\'eodory number of
$V^{1,2}$ is at most $5$.

Given any number $N$ of tensors $(\chi_1,\ldots,\chi_N)$ we may
define the convex hull $V^{1,\ldots,N}$. In
Subsection~\ref{S:dimension} we prove that $\dim V^{1,\ldots,N}$ is
$5$ times the dimension of the span of $(\chi_1,\ldots,\chi_N)$. 
In the text, we will omit the superscripts from our notation for the
convex hull.

Magnetic susceptibility tensors are $3\times 3$ symmetric
trace zero matrices and form a 5-dimen\-sion\-al irreducible representation of the
group $SO(3)$ of rotations in $\R^3$.
More generally, one could study the convex hulls of orbits of compact groups.
We were surprised to find that very little is known about such convex bodies,
particularly their Carath\'eodory numbers and facets.
We hope that our work will stimulate a more thorough study of convex
hulls of orbits of compact groups.

In Section~\ref{S:motivation}, we describe the motivation for this work
from protein structure.
In Section~\ref{S:group} we discuss group actions and in Section~\ref{S:convex}
convex hulls of orbits.
In Section~\ref{S:One} we complete the analysis of~\cite{GLS05} in the case
of one metal ion.
In Section~\ref{S:Two} we analyze the case of two metal ions and deduce
Theorem~\ref{T:coaxial}.

%
%
\section{Application to protein structure}\label{S:motivation}

Proteins are large biological molecules synthesized by living
organisms. The genome of an organism contains the chemical formulae
for its proteins. Currently, hundreds of organisms (including man)
have had their genomes mapped, and such chemical formulae are
readily available. An important step towards inferring the
biological function of a protein from its chemical formula is to
determine its 3-dimensional structure, or its \Blue{{\it fold}}.

About one third of all proteins 
incorporate metal ions into their structures. The fold of these
proteins may be inferred from nuclear magnetic resonance, which can
measure the interactions between paramagnetic metal ions and dipoles
within the protein. The main quantities that can be measured are the
pseudo contact shifts (PCS) \cite{BBBCGLT} and the residual dipolar
coupling (RDC) \cite{TFKP}.
In this paper we only deal with the RDC.

The residual dipolar coupling between a paramagnetic ion and a
dipole formed by atoms $a$ and $b$ within the protein depends upon
the vector displacement $r$ from the atom $a$ to the atom $b$ and
the \Blue{{\em magnetic susceptibility tensor} $\chi$} of the metal
ion, which is a $3\times 3$ symmetric matrix. The RDC interaction
has the following vector formula
 \begin{equation} \label{eq1}
   \delta\ :=\  \frac{C}{\|r\|^5}\, r^T \chi r\ -\
       \frac{C}{3\|r\|^3}\mbox{Trace}(\chi)\,.
 \end{equation}
Here, $C$ is a constant and $\|r\|$ is the length of the vector $r$.
This only depends upon the relative orientation of the dipole and metal
ion, and so the RDC data may be used to infer this relative orientation.

Writing $\chi=\chi_0 + \frac{1}{3}\mbox{Trace}(\chi)I_3$, where $I_3$ is the
$3\times 3$ identity matrix and $\chi_0$ is the trace-free or
\Blue{{\it anisotropic}} part of $\chi$, this formula becomes
\[
    \delta\ =\  \frac{C}{\|r\|^5}\, r^T \chi_0 r\,.
\]
We assume henceforth that $\chi=\chi_0$ is anisotropic.

\begin{figure}[htb]
\begin{picture}(326,135)(-50,0)
   \put(0,0){\includegraphics[height=4.5cm]{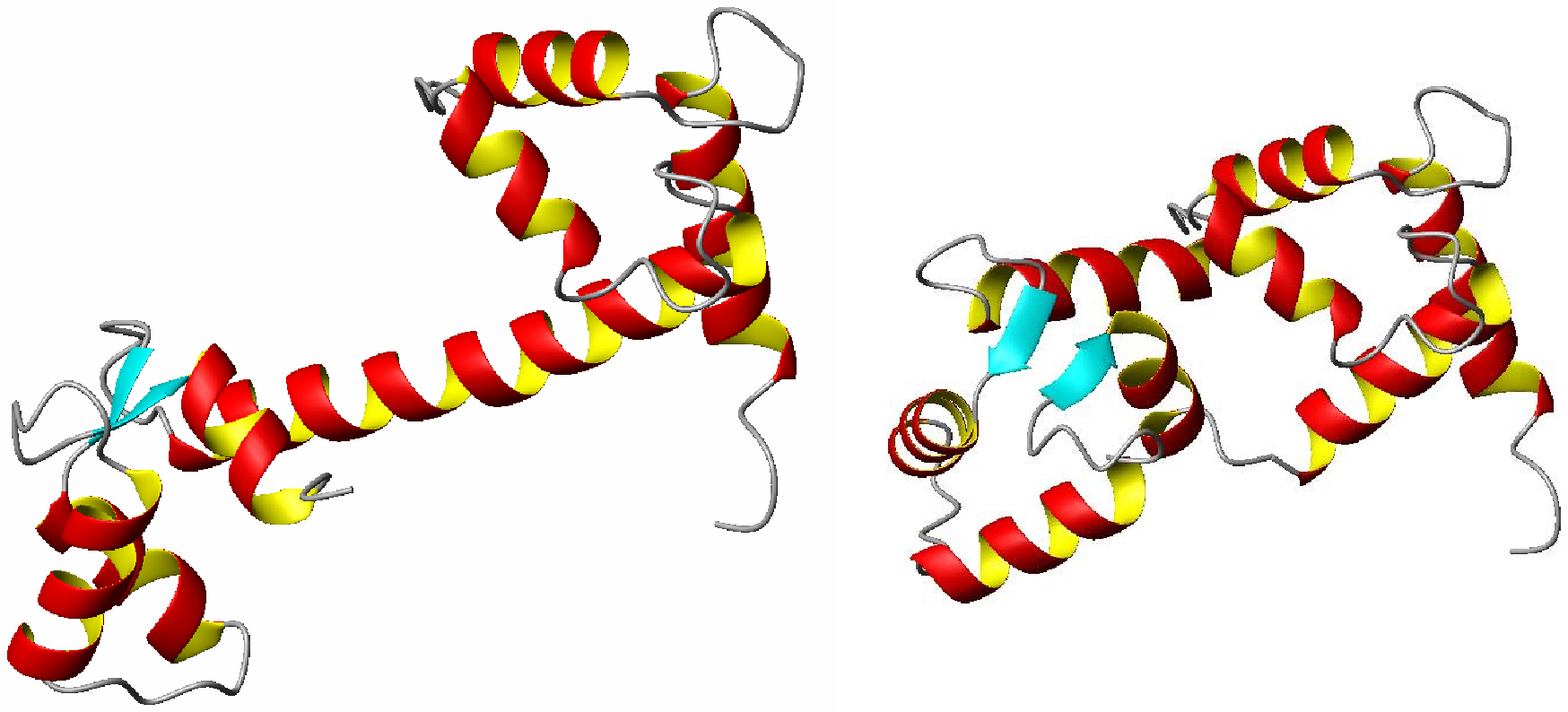}}
   \put(120,0){linker zone}
   \put(85, 3){\vector(0,1){50}}  \put(85, 3){\line(1,0){30}}
   \put(180, 3){\line(1,0){30}} \put(210, 3){\vector(0,1){40}}

   \put(-50,85){C-terminal} \put(-17,82){\vector(2,-1){30}}
   \put(153,132){N-terminal}
   \put(165,128){\vector(-2,-1){20}}  \put(200,128){\vector(1,-1){20}}

   \put(278,77){\begin{minipage}{2in}ion\\binding\\site\end{minipage}}
   \put(275,77){\vector(-1,0){35}}
 \end{picture}
 \caption{Two orientations of calmodulin}\label{calmod}
\end{figure}

The fold of the protein is usually unique, in the sense that small
variations of the shape are allowed. There are proteins which
however exhibit large variations of shape under particular
conditions. A widely-studied example is calmodulin, which has two
rigid domains, called the \Blue{{\it N-terminal}} and \Blue{{\it
C-terminal}} domains, connected by a short flexible linker. The N-
and C-terminal domains are assumed to be rigid bodies with known
structures. Figure~\ref{calmod}, obtained with Molmol \cite{KBW},
shows calmodulin in two different orientations.

The binding site of the metal ion in calmodulin belongs to the
N-terminal domain. The measured RDC of pairs of atoms belonging to
the N-terminal domain can be used to obtain a good estimate of
$\chi$. The measured RDC of pairs of atoms belonging to the
C-terminal domain can be used to study the relative orientation of
the two domains.

Let us model the relative orientation of the N- and C-terminal domains with a
rotation $R$.
Then there is an unknown probability measure $p$ on the set $SO(3)$ of rotations
such that the \Blue{{\em mean RDC} $\overline{\delta}$} of the pair of atoms $a,b$ in the
C-terminal domain is given by
 \begin{equation}\label{calcdata}
   \overline{\delta}\  =\  \frac{C}{\|r\|^5}\int_{SO(3)} {(Rr)^T \chi (Rr)\,dp(R)}\,
    \ =\ \frac{C}{\|r\|^5}r^T \chibar r\,,
 \end{equation}
where the \Blue{\textit{mean magnetic susceptibility tensor} $\chibar$} is
 \begin{equation}\label{genmtensor}
   \chibar\ =\ \int_{SO_3} R^T\chi R\,dp(R)\,.
 \end{equation}
This tensor $\chibar$ can be estimated from the RDC of several dipole pairs in
the C-terminal domain.
The experimental measures show that in terms
of difference of eigenvalues, $\chibar$ is between $5$ and $20$ times smaller
than $\chi$~\cite{BDGKLPPPZ}.
This indicates that $p$ is not a point mass, that is, the C-terminal domain
moves with respect to the N-terminal domain.

The availability of $N$ distinct mean
susceptibility tensors $\chibar_k$ with respect to different metal
ions $k=1,\ldots,N$ increases the information about $p$ for $N$ up
to $5$, see for instance \cite{MPPGB} and~\cite[Theorem~3.2]{LLPS}.
However, even the exact knowledge of $5$ mean tensors $\chibar_k$
(i.e. $25$ real numbers) does not allow the exact reconstruction of
the probability measure $p$.

An approach to extract information from the mean tensors is to
define $p_{\max}(R)$ as the maximal fraction of time that the
C-terminal can stay in a particular orientation $R$, yet still
produce the measured mean tensors. Orientations with a large
$p_{\max}$ agree with with what are thought to be the most favored
orientations of the C-terminal~\cite{LLPS}.

The calculation of $p_{\max}$ can be performed geometrically if only
RDC is considered \cite{LLPS}. In the combined PCS+RDC case more
information is added \cite{BGLPPSY}, however the calculation of
$p_{\max}$ can only be performed numerically. For the convergence
and efficiency of the algorithm, the minimal number of orientations
needed to reconstruct any admissible set of mean tensors $\chibar_k$
should be used. Experience suggests that adding the PCS data does
not increase the actual number of orientations needed
\cite{BGLPPSY}, so the Carath\'eodory number for the RDC case may be
used as a basis for the numerical minimization.

%
%
\section{Groups actions and anisotropic tensors}\label{S:group}
%

We first recall some basics about representations of compact groups, then consider
the action of the group $SO(3)$ of rotations in $\R^3$ on the 5-dimensional space of
anisotropic $3\times 3$ symmetric tensors, and finally investigate the coaxial subgroups of
$SO(3)$.

%
\subsection{Representations of compact groups}
%
This material may be found in the book of
Br\"ocker and tom Dieck~\cite[Ch. II]{BtD}.
Let $G$ be a compact group, such as $SO(3)$.
A \Blue{{\it representation}} of $G$ is a finite-dimensional vector space $W$ on
which $G$ acts by linear transformations.
That is, we have a group homomorphism $\rho\colon G\to GL(W)$, where
$GL(W)$ is the group of invertible linear transformations on $W$.
For $g\in G$ and $w\in W$, write $g.w$ for $\rho(g)(w)$.

A representation $W$ of $G$ is \Blue{{\it irreducible}} if its only
$G$-invariant subspaces are $\{0\}$ and $W$.
Every representation of $G$ decomposes as the direct sum of irreducible
representations which is unique in the following way.
Given a representation $W$ of $G$ and a positive integer $l$,
let $W^l$ be the $l$-fold direct sum of $W$,
\[
   W^l\ =\ \underbrace{W\oplus W\oplus\dotsb \oplus W}_l\,.
\]
Elements $g$ of $G$ act diagonally on elements $w=(w_1,\dotsc,w_l)$
of $W^l$, $g.w=(g.w_1, \dotsc,g.w_l)$. Suppose that $W_1, W_2,
\dotsc$ is the list of irreducible representations of $G$. If $U$ is
a representation of $G$ then there exist unique integers
$l_1,l_2,\dotsc$ such that
\[
   U\ \simeq\  W_1^{l_1} \oplus W_2^{l_2}\oplus W_3^{l_3}\oplus \dotsb\ ,
\]
as representations of $G$.
If $U_i$ is the subrepresentation of $U$ mapped to the summand $W_i^{l_i}$ under this
isomorphism, then $U_i$ does not depend on any choices and is called the
\Blue{{\it isotypical component}} of $U$ corresponding to $W_i$.
If $l_i>0$, then we say that $U$ \Blue{{\it contains}} $W_i$.
Furthermore, if $U'\subset U$  is a subrepresentation, then the $i$th isotypical component
of $U'$ is $U'\cap U_i$, which is also the image of $U'$ under the projection to
$U_i$.

\Blue{{\it Haar measure}} is a $G$-invariant measure \Blue{$\mu$} on $G$ with
$1=\int_Gd\mu(g)$.
Given a linear function $L\colon W\to \R$, where $W$ is a representation of $G$, we may
average $L$ over orbits of $G$ to get a new linear function $L'$, defined by
\[
   \Blue{L'}(x)\ :=\ \int_G L(g.x)d\mu(g)\,.
\]
Since $L'$ is constant on orbits of $G$, it is $G$-invariant.
This association $L\mapsto L'$ is called the \Blue{{\it Reynolds operator}}.
It is an important tool for analyzing $G$-representations.

Another key tool is Schur's lemma.
A linear map $\varphi\colon W\to U$ between representations of $G$ is a
\Blue{{\it $G$-map}} if for all $w\in W$ and $g\in G$, we have
$g.\varphi(w)=\varphi(g.w)$.
Let \Blue{$\Hom_G(W,U)$} be the space of $G$-maps.
A \Blue{{\it division algebra}} is a finite-dimensional associative algebra in
which every non-zero element is invertible.
\medskip

\noindent{\bf Schur's Lemma.}
{\it
  If\/ $W\not\simeq U$ are irreducible representations of $G$, then $\Hom_G(W,U)=0$
  and $\Hom_G(W,W)$ is a division algebra which contains $\R$.
}\medskip

\noindent{\it Proof.}
Let $\varphi\colon W\to U$ be a $G$-map.
Then the kernel of $\varphi$ is a subrepresentation of $W$ and so it is
either $0$ or $W$, and the image of $\varphi$, which is a
subrepresentation of $U$, is either 0 or $U$.
Examining the possibilities leads to the conclusions.
\QED\smallskip

There are exactly 3 division algebras which contain $\R$:
The real numbers $\R$, the complex numbers $\C$, and the quaternions
$\H$.
An irreducible representation $W$ of $G$ has \Blue{{\it real}},
   \Blue{{\it complex}}, or \Blue{{\it quaternionic type}}, depending on
$\Blue{\mbox{End}_G(W)}:=\Hom_G(W,W)$.

\begin{example}
 Consider the group $SO(2)$ of rotations of $\R^2$,
\[
   SO(2)\ =\ \left\{ \Blue{R_\theta}\ :=\ {\textstyle   \left(
  \begin{array}{rr}\cos\theta&-\sin\theta\\\sin\theta&\cos\theta\end{array}
   \right)}\ \mid\ \theta\in [0,2\pi)\right\}\ .
\]
$\mbox{End}_{SO(2)}(\R^2)$ consists of those $2\times 2$ matrices $M$ such that
$MR_\theta=R_\theta M$, and so
\[
   \mbox{End}_{SO(2)}(\R^2)\ =\
    \R{\textstyle \left(\begin{matrix}1&0\\0&1\end{matrix}\right)}
    \ +\ \R\left(\!\begin{array}{rr}-1&0\\0&1\end{array}\right)\ .
\]
This is isomorphic to $\C$ (we send $\refl$ to $\sqrt{-1}$),
so this representation of $SO(2)$ has complex type.
It is the \Blue{{\it defining representation} $U_1$} of $SO(2)$.
If we identify $\R^2$ with $\C$ and $SO(2)$ with the circle group
$\Blue{S^1}:=\{e^{i\theta}\mid 0\leq\theta<2\pi\}$, then the action
is scalar multiplication by elements of $S^1$.
For any positive integer $k>0$, let \Blue{$U_k$} be the representation of $S^1$ on $\C$
(identified with $\R^2$) where $z\in S^1$ acts as multiplication by $z^k$.
These all have complex type.
\end{example}

\begin{example}
 The orthogonal group \Blue{$O(2)$} contains $SO(2)$ as well
 as the cosets of reflections
\[
  SO(2)\cdot \left(\begin{array}{rr}0&-1\\1&0\end{array}\right)\ =\
  \left\{\ R_\theta\cdot
   \left(\!\begin{array}{rr}-1&0\\0&1\end{array}\right)\ \mid\
   \theta\in[0,2\pi)\right\}\ .
\]
 The defining representation \Blue{$U_1$} of $O(2)$ on $\R^2$ has real
 type, as
\[
   \left(\!\begin{array}{rr}-1&0\\0&1\end{array}\right)\cdot
   \left(\begin{array}{rr}0&-1\\1&0\end{array}\right)\ =\
   \left(\begin{matrix}0&1\\1&0\end{matrix}\right)\ \neq\
   \left(\begin{array}{rr}0&-1\\-1&0\end{array}\right)\ =\
   \left(\begin{array}{rr}0&-1\\1&0\end{array}\right)\cdot
   \left(\!\begin{array}{rr}-1&0\\0&1\end{array}\right)\ ,
\]
  and so  $\mbox{End}_{O(2)}(\R^2)=\R\cdot\iden\simeq \R$.

In the \Blue{{\it trivial representation}}  $\Blue{U_0}=\R$ of $O(2)$, elements act as
multiplication by 1.
For a positive integer $k$, define the map
$\varphi_k\colon O(2)\to O(2)\subset GL(2,\R)$ by
$\varphi_k(R_\theta)=R_{k\theta}$ and $\varphi_k\refl=\refl$.
This defines the representation \Blue{$U_k$} of $O(2)$, which has real type.
Restricting to $SO(2)$ gives its representation $U_k$ of complex
type.
\end{example}

%
\subsection{Rotations of anisotropic tensors}
%

Let $e_1=(1,0,0)^T$, $e_2=(0,1,0)^T$, and $e_3=(0,0,1)^T$ be the standard basis
of $\R^3$.

The \Blue{{\it special orthogonal group $SO(3)$}} is the group of rotations in $\R^3$.
It consists of $3\times 3$ real orthogonal matrices
with determinant 1,
 \[
   SO(3)\ :=\ \left\{R\in \R^{3\times 3}\mid RR^T=1\quad\mbox{and}\quad \det R=1\right\}\,.
 \]
Let $R\in SO(3)$, and let $T_R SO(3)$ be the tangent space to
$SO(3)$ at the matrix $R$. Let $I$ be the identity matrix, then $T_I
SO(3)$ is the space $\mathfrak{so}_3$ of skew symmetric $3\times 3$
matrices, which is the Lie algebra of $SO(3)$. That is,
\[
    T_I SO(3)\ =\ I + \mathfrak{so}_3\,.
\]
Elements $R\in SO(3)$ act on $3\times 3$ symmetric matrices
(tensors) $\chi$ by conjugation, $R.\chi := R\chi
R^T$\fauxfootnote{This left action (if $R,S\in SO(3)$, then
  $R(S.\chi)=RS.\chi$) is equivalent to the implied action in~\eqref{genmtensor}.}.
This preserves the trace of $\chi$, and so $SO(3)$ acts on the
space $W$ of \Blue{{\em anisotropic}} (trace-zero) tensors,  a 5-dimensional irreducible
real representation.
We introduce some useful coordinates for $W$.
A point $(v,w,x,y,z)\in\R^5$ corresponds to 
 \begin{equation}\label{Eq:Coords}
  \chi(v,w,x,y,z)\ :=\
  \left(\begin{array}{crr}v&0&0\\0&-\tfrac{v}{2}&0\\0&0&-\tfrac{v}{2}\end{array}\right)
   \ + \
  \left(\begin{matrix}0&w&x\\w&0&0\\x&0&0\end{matrix}\right)
   \ + \
  \left(\begin{array}{ccr}0&0&0\\0&y&z\\0&z&-y\end{array}\right)\ .
 \end{equation}
%
Observe that $e_1$ is an eigenvector for $\chi(v,w,x,y,z)$ if and only if $w=x=0$,
$e_2$ is an eigenvector if and only if $w=z=0$,  and
$e_3$ is an eigenvector if and only if $x=z=0$.

%
\subsection{Coaxial subgroups}\label{S:coaxial}
%

The \Blue{{\it coaxial subgroup $Q_e$}} is the set of rotations
fixing a line in $\R^3$ with direction $e$. It is isomorphic to the
orthogonal group $O(2)$. Its identity component $Q_e^+$ is
isomorphic to $SO(2)$ and consists of rotations about the axis $e$,
while the other component $Q_e^-$ consists of reflections in axes
orthogonal to $e$. For example, let
 \begin{equation}\label{Eq:coaxial_members}
   R_{e_1,\theta}\ :=\ \left(\begin{matrix}1&0&0\\0&\cos\theta&-\sin\theta
         \\0&\sin\theta&\quad \cos\theta\end{matrix}\right)
   \qquad\mbox{and}\qquad
   R_{e_3,\pi}\ :=\ \left(\begin{array}{rrc}-1&0&0
              \\0&-1&0\\0&0&1\end{array}\right)\ .
 \end{equation}
If we fix $f$ perpendicular to $e$ and let $\theta$ run over all
angles, then $R_{e,\theta}$ and $R_{e,\theta}R_{f,\pi}$ give all
elements of $Q_e$.

We consider the action of a coaxial subgroup $Q_e$ on $W$.
For this, suppose that $e=e_1$ and let $R_{e_1,\theta}$ act on
$\chi(v,w,x,y,z)$.
This gives the tensor $\chi(v',w',x',y',z')$, where
 \begin{eqnarray*}
   v'&=&  v\ ,\\
   \left(\begin{matrix}w'\\x'\end{matrix}\right)&=&
    \left(\begin{array}{rr}\cos\theta&-\sin\theta\\
          \sin\theta&\cos\theta\end{array}\right)
         \left(\begin{matrix}w\\x\end{matrix}\right)\quad=\quad
       R_\theta \ \left(\begin{matrix}w\\x\end{matrix}\right)\,,\quad\mbox{and}\\
   \left(\begin{matrix}y'\\z'\end{matrix}\right)&=&
    \left(\begin{array}{rr}\cos2\theta&-\sin2\theta\\
          \sin2\theta&\cos2\theta\end{array}\right)
         \left(\begin{matrix}y\\z\end{matrix}\right)\quad=\quad
       R_{2\theta}\left(\begin{matrix}y\\z\end{matrix}\right)\,.
 \end{eqnarray*}
Thus $R_{e_1,\theta}$ acts trivially
on the coordinate $v$, by rotation through the angle $\theta$ on the
vector $(w,x)^T$, and by rotation through $2\theta$ on the vector
$(y,z)^T$.
Note that $R_{e_3,\pi}$ sends $\chi(v,w,x,y,z)$ to $\chi(v,\,w,-x,\,y,-z)$.
Thus, if we restrict the action of $SO(3)$ on $W$ to
$Q_{e_1}\simeq O(2)$, then it
decomposes as a sum of irreducible representations
 \begin{equation}\label{Eq:Decomposition}
   W\ =\ U_0\ \oplus\ U_1\ \oplus\ U_2\,.
 \end{equation}
This decomposition corresponds to the coordinates~\eqref{Eq:Coords}.
Projection to the trivial submodule $U_0=\R$ is, up to a scalar multiple, the
unique $Q_e$-invariant linear function $L\colon W\to \R$.

%
%
\section{Convex hulls of orbits}\label{S:convex}
%

Let $\calO$ be an orbit of a compact group $G$ in a representation
$W$ of $G$. The \Blue{{\it convex hull} $V$} of $\calO$ is all
points of $W$ which are convex combinations of elements of $\calO$,
 \[
   \lambda_1v_1+\lambda_2v_2+\dotsb+\lambda_nv_n\,,
 \]
where $v_1,\dotsc,v_n\in \calO$, and the non-negative numbers
$\lambda_i$ have sum 1. The set $V$ is a compact convex set, hence a
convex body.

%
\subsection{Faces and Carath\'eodory number of $V$}
%

%
A \Blue{{\it face}} $F$ of $V$ is the subset of its boundary
where some linear function $L$ achieves its maximum on $V$,
\[
   F\ :=\ \left\{v\in V\mid L(v)\geq L(u)
     \mbox{ for all }u\in V\right\}\,.
\]
We say that $L$ \Blue{{\it supports}} $F$ and also that the
hyperplane $L(x)=L(F)$  \Blue{{\it supports}} $F$. (Here, $L(F)$ is
a constant.) The tangent spaces to $\calO$ of its points lying in
$F$ are contained in any hyperplane supporting $F$. Such a tangent
space \Blue{$T_v\calO$} at the point $v$ is
\[
   v\ +\ \mathfrak{g}.v\,,
\]
where the action of the Lie algebra $\mathfrak{g}$ is the derivative
of the action of $G$. The face $F$ is \Blue{{\it proper}} if $F\neq
V$. When $V$ is full-dimensional so that $\dim V=\dim W$, this is
equivalent to $L\neq 0$. A \Blue{{\it facet}} is a maximal face and
a \Blue{{\it vertex}} is a minimal face. Vertices are not a convex
combination of other points of $V$.

\begin{lemma}
  The vertices of\/ $V$ are exactly the points of\/ $\calO$.
\end{lemma}

\noindent{\it Proof.}
 The vertices of $V$ are a subset of $\calO:=G.x$.
 Let $g.x\in\calO$ and suppose that it is a convex combination of vertices,
\[
    g.x\ =\ \lambda_1 g_1.x + \lambda_2 g_2.x + \dotsb +\lambda_n g_n.x \,.
\]
 Multiplying by $g_1g^{-1}$ expresses the vertex $g_1.x$ as a convex combination of points
 of $\calO$.
 Thus $n=1$ and $g.x=g_1.x$ is a vertex.
\QED\smallskip

When $V$ has dimension $d$,  Carath\'eodory's Theorem~\cite{Ca1911},
see e.g. \cite[Theorem~1.1.4]{Schneider}, implies that any point $x$
of $V$ is a convex combination of at most $d{+}1$ vertices. The
\Blue{{\it Carath\'eodory number}} of $V$ is the minimum number $n$
such that any point $x\in V$ is a convex combination of at most $n$
vertices. For example, a ball in $\R^d$ has Carath\'eodory number
$2$, while a $d$-simplex has Carath\'eodory number $d+1$.
Fenchel~\cite{Fe1929}, see e.g. \cite[Theorem~1.4]{Reay}, showed
that the Carath\'eodory number is at most $d$ when the set of
vertices is connected. More useful for us is a recursive bound,
which is immediate from the observation that any point of $V$ is the
convex combination of any vertex and some boundary point.

\begin{lemma}\label{L:Cartheodory}
 The Carath\'eodory number of a convex body $V$ is at most one more than the
 maximal Carath\'eodory number of its facets.
\end{lemma}

Suppose that $S$ is a closed (hence compact) subgroup of $G$ which stabilizes
a face $F$ of $V$, that is $s.F=F$ for all $s\in S$.

\begin{lemma}\label{Lem:trivial_support}
  When $F$ is proper, there is a non-zero $S$-invariant linear function on $W$, and $W$
  contains the trivial representation of $S$.
\end{lemma}

\noindent{\it Proof.}
 Let $L\colon W\to \R$ be any linear function supporting $F$ with $L(F)=\ell$.
 Let $L'$ be the image of $L$ under the Reynolds operator for $S$.
 For $u\in F$, we have $L'(u)=\ell$, as $F$ is $S$-stable and $L'(u)$ is the average of
 $L$ over the orbit of $S$ through $u$.

 Suppose that $w\in V \setminus F$.
 Then $L(w)<L(u)=\ell$ and $L(s.w)$ is bounded away from $\ell$ as
 $S.w$ is compact and disjoint from $F$.
 In particular, this implies that $L'(w)<L'(u)$, which shows that the $S$-invariant linear
 function $L'$ supports $F$ and that $L'\neq 0$.
 Such an $S$-invariant linear function must factor through the trivial isotypical
 component of $W$ as a representation of $S$.
 This completes the proof.
\QED

%
\subsection{The dimension of $V$}\label{S:dimension}
%
If $w=(w_1,\dotsc,w_l)\in W^l$, then we write $d(w)$ for the
dimension of the linear span of the components $w_1,\dotsc,w_l$ of
$w$ in $W$.

\begin{lemma}\label{Lem:orbit-Span}
 Suppose that $W$ is an irreducible representation of a group $G$ having real type.
 If $w\in W^l$, then the linear span $U$ of the orbit $G.w$ in $W^l$ is
 isomorphic to $W^{d(w)}$.
\end{lemma}

\noindent{\it Proof.}
 Write $k:=d(w)$.
 We may assume that
 $w_1,\dotsc,w_k$ form a basis for the linear span of $w_1,\dotsc,w_l$.
 Let $A=(\alpha_{ij})\in\Mat_{l\times k}(\R)$ be the matrix which writes
 the components of $w$ in terms of this basis,
 $w_{i} = \sum_{j=1}^k \alpha_{i,j}w_j$ for $i=1,\dotsc,l$.
 For each $i=1,\dotsc,l$, let
 $\varphi_i\colon U\to W$ be the projection to the $i$th coordinate.
 Since $\varphi_i(g.w)=g.w_i$, we have
 \begin{equation}\label{Eq:phi_sum}
  \varphi_{i}\ =\ \sum_{j=1}^k \alpha_{i,j}\varphi_j
   \qquad\mbox{for}\quad i=1,\dotsc,l\,.
 \end{equation}
 This matrix $A$ defines a $G$-map  $A\colon W^k\to W^l$ by
 \begin{equation}\label{Eq:A_map}
  A\ \colon\ (w_1,\dotsc,w_k)\ \longmapsto\
    \left(\sum_j \alpha_{1,j}w_j,\,\sum_j \alpha_{2,j}w_j,\,
           \dotsc,\,\sum_j \alpha_{l,j}w_j\right)\ .
 \end{equation}
 Composing the map $\psi:=(\varphi_1,\dotsc,\varphi_k)\colon U\to W^k$ with
 $A\colon W^k\to W^l$ gives the identity map on $U$:
 By~\eqref{Eq:phi_sum}, for $w\in U$, we have
\[
   w\ =\ (\varphi_1(w),\varphi_2(w),\dotsc,\varphi_l(w))\,.
\]

We show that the map~\eqref{Eq:A_map} is injective and thus $\psi$ is an
isomorphism.
A linear map $L\colon W\to \R$ induces maps $L^k\colon W^k\to\R^k$ and
$L^l\colon W^l\to\R^l$, which commute with $A$.
If $0\neq w\in W^k$, then there is some linear map $L\colon W\to \R$
with $L^k(w)\neq 0$.
Since $A$ has full rank $k$, $A(L^k(w))\neq 0$.
But this implies that $A(w)\neq 0$, as $A(L^k(w))=L^l(A(w))$.
\QED

\begin{lemma}\label{L:Orbit_dimension}
  Suppose that $W=W_1^{l_1}\oplus\dotsb\oplus W_m^{l_m}$ is the decomposition of
  a representation $W$ of $G$ into isotypical pieces, each of which has real type.
  Let  $w=(w_1,\dotsc,w_m)\in W$, where  $w_i$ is the component of $w$ in $W_i^{l_i}$.
  Then the dimension of the convex hull $V$ of the orbit $G.w$ is
\[
    \sum_{W_i\neq \R}  d(w_i)\cdot \dim W_i\,.
\]
  If $W$ does not contain the trivial representation, then $0$ lies in $V$.
\end{lemma}

\noindent{\it Proof.}
 If $W$ contains the trivial representation, assume that it is $W_1$.
 Since $g.(w-w_1)=g.w-w_1$, we see that the orbits
 $G.w$ and $G.(w-w_1)$ are isomorphic, and the same it true for their convex hulls.
 Thus it it no loss to suppose that $w_1=0$, which is equivalent to assuming that $W$ does
 not contain the trivial representation.

 The linear span $U$ of the orbit $G.w$ is the direct sum of its projections to the
 isotypical components $W_i^{l_i}$ of $W$.
 Each projection is the linear span of $G.w_i$, which by Lemma~\ref{Lem:orbit-Span} is
 isomorphic to $W_i^{d(w_i)}$.
 Thus $U$ has dimension  $\sum_i d(w_i)\dim W_i$.
 We may replace $W$ by this linear span, and therefore assume that the orbit $G.w$ spans
 $W$.

 The convex hull of $G.w$ lacks full dimension only if it lies in some hyperplane $H$ not
 containing the origin.
 Suppose that this is not the case and let $B$ be the convex hull of $G.w$ and the
 origin.
 Then $V$ is a proper $G$-stable face of $B$ and so by Lemma~\ref{Lem:trivial_support}
 $W$ contains the trivial representation, which is a contradiction.

 If $0\not\in V$, then there is some linear function $L$ which is bounded above 0 on $V$.
 But then the image $L'$ of $L$ under the Reynolds operator is non-zero on $V$.
 This implies that $L'\neq 0$, and so $W$ contains the trivial representation,
 a contradiction.
\QED

\begin{example}\label{Ex:SO(2)_Orbits}
 Lemmas~\ref{Lem:orbit-Span} and~\ref{L:Orbit_dimension}
 do not hold if the representation $W$ has complex type.
 For example, let $G=SO(2)$ and $W=U_k^l$ with $k,l\geq 1$.
 Identifying $W$ with $\C^l$ and $SO(2)$ with the circle group $S^1$, elements $z\in S^1$
 act on $\C^l$ as scalar multiplication by $z^k$.
 Thus the linear span $SO(2).w$ for $w\in U_k^l$ is
 a complex line, and therefore has real dimension 2, and not 4 as
 Lemma~\ref{Lem:orbit-Span} predicts for general $w\in W$ when $l,k\geq 1$.

 In particular, if $W=U_1\oplus U_2^2$, and $w\in W$ is general, then
 the linear span of $G.w$ has complex dimension 2 and thus real dimension 4.
\end{example}

%
\section{One metal ion}\label{S:One}
%

Let $W$ be the space of symmetric $3\times 3$ anisotropic tensors, a
5-dimensional irreducible representation of $SO(3)$ of real type.
For each unit vector $e\in\R^3$, there is a linear function
 \begin{equation}\label{E:L_e}
   \Blue{L_e}\ \colon\ W\ni\chi \ \longmapsto\
   \langle e, \chi e\rangle\ =\ e^T\chi e\in\R\,.
 \end{equation}
If $e$ is an eigenvector of $\chi$, then $L_e(\chi)$ is its
eigenvalue. In general, $L_e(\chi)$ lies between the maximum and
minimum eigenvalues of $\chi$. Note that $L_e$ is $Q_e$-invariant.
By the decomposition~\eqref{Eq:Decomposition} of $W$ into
irreducible $Q_e$ representations, any $Q_e$-invariant linear
function is a scalar multiple of $L_e$.

As a matrix, a tensor in $W$ has real eigenvalues and its
eigenvectors form an orthonormal basis for $\R^3$. Fix a non-zero
anisotropic tensor  $\chi\in W$ with maximum eigenvalue $M>0$ and
minimum eigenvalue $m<0$. The intermediate eigenvalue of $\chi$ is
$-M-m$, and we have $-\frac{M}{2}\geq m\geq -2M$. The orbit $\calO$
of $\chi$ under $SO(3)$ consists of the anisotropic tensors with
maximal eigenvalue $M$ and minimal eigenvalue $m$.
It is a manifold whose dimension we determine.

\begin{proposition}\label{Prop:orbit_dim}
 The orbit $\calO_\chi$ is $3$-dimensional unless $\chi$ has an eigenvalue of
 multiplicity $2$ and then it is two-dimensional.
\end{proposition}

\noindent{\it Proof.}
The dimension of $\calO$ is equal to the dimension of any of its tangent
spaces.
Since
\[
   \calO\ =\ \left\{ R\chi R^T\mid R\in SO(3)\right\}\,,
\]
the tangent space $T_\xi\calO$ at a point $\xi\in\calO$ is the affine space
 \begin{equation}\label{Eq:TS}
   \xi\ +\ \left\{ r\xi + \xi r^T\mid r\in\mathfrak{so}_3\right\}\,.
 \end{equation}
Indeed, consider the action in $\xi$ of an element $I+r$ of the tangent space
$T_ISO(3)=I+\mathfrak{so}_3$:
\[
   (I+r)\xi(I+r)^T\ =\ \xi + r\xi + \xi r^T + r\xi r^T  \ .
\]
Discarding the term which is quadratic in $\mathfrak{so}_3$ gives~\eqref{Eq:TS}.

It suffices to determine the tangent space to $\calO$ at the point $\chi$.
We may suppose that $\chi$ is diagonal, and let $r$ be a general
element of $\mathfrak{so}_3$,
\[
   \chi\ =\ \left(\begin{matrix}M&0&0\\0&-M-m&0\\0&0&m\end{matrix}\right)
   \qquad\mbox{and}\qquad
   r\ =\ \left(\begin{matrix}0&-a&-b\\a&0&-c\\b&c&0\end{matrix}\right)\ ,
\]
where $a,b,c\in\R$.
Let $\alpha:=M+\frac{m}{2}$, $\beta:=M-m$, and $\gamma:=-m-\frac{M}{2}$.
Then $\beta>0$ and $\alpha,\gamma\geq 0$ with $\alpha=0$ only when the eigenvalue $m$
has multiplicity 2 and $\gamma=0$ only when the eigenvalue $M$ has multiplicity 2.
We see that $T_\chi\calO$ is the affine subspace of $W$
 \begin{equation}\label{Eq:TanSpace_k=1}
   \xi
    \ +\
    \left(\begin{matrix}0&2a\alpha&b\beta\\
                     2a\alpha&0&2c\gamma\\b\beta&2c\gamma&0\end{matrix}\right)\ ,
 \end{equation}
where $a,b,c\in\R$.
This is 3-dimensional unless either $\alpha=0$ or $\gamma=0$.
\QED\smallskip

Let $V$ be the convex hull of the orbit $\calO$ of $\chi\in W$.
By Lemma~\ref{Lem:orbit-Span}, this is a 5-dimensional convex body.

\begin{lemma}\label{Lem:Le-bounds}
  If\/ $\chibar\in V$ and $e\in\R^3$ is a unit vector, then we have
 \begin{equation}\label{E:Le_ineq_orbit}
   m\ \leq\ L_e(\chibar)\ \leq\ M\,.
 \end{equation}
\end{lemma}
In fact, $V$ is the set of symmetric anisotropic tensors
satisfying~\eqref{E:Le_ineq_orbit}~\cite[Theorem~3.3]{GLS05}.\smallskip

\noindent{\it Proof.}
 $L_e(\chibar)$ lies between the maximum and minimum eigenvalues of $\chibar$.
 Thus the inequality~\eqref{E:Le_ineq_orbit} holds for $\chibar$ in the orbit of $\chi$.
 Since a general element of $V$ is a convex combination of tensors in the orbit
 of $\chi$, we deduce~\eqref{E:Le_ineq_orbit}.
\QED\smallskip

A \Blue{{\it coaxial face}} of $V$ is a face which is stabilized by some
coaxial subgroup $Q_e$.
By Lemma~\ref{Lem:trivial_support}, a coaxial face stabilized by
$Q_e$ is supported by a non-trivial $Q_e$-invariant linear function.
As we noted earlier, this linear function is necessarily a scalar multiple of $L_e$.
By the inequality~\eqref{E:Le_ineq_orbit}, there are two possibilities for
such a coaxial face,
 \begin{equation}\label{Eq:Coax_Faces}
  \begin{array}{rcl}
   \Blue{F_e^M}&:=& \left\{\chibar\in V\mid L_e(\chibar)=M\right\}\quad\mbox{and}\\
   \Blue{F_e^m}&:=& \left\{\chibar\in V\mid L_e(\chibar)=m\right\}\,.
    \rule{0pt}{15pt}
  \end{array}
 \end{equation}
The coaxial face $F^M_e$ consists of tensors $\chibar\in V$ having $e$
as an eigenvector with eigenvalue $M$ and tensors in $F^m_e$ have $e$ as an
eigenvector with eigenvalue $m$.

We now describe the facets of $V$ and determine its Carath\'eodory number.
As in the proof of Proposition~\ref{Prop:orbit_dim}, set
$\Blue{\alpha}:=M+\frac{m}{2}\geq 0$ and $\Blue{\gamma}:=-m-\frac{M}{2}\geq 0$.

\begin{lemma}\label{L:Boundary_one}
 The boundary of\/ $V$ is the union of coaxial faces~\eqref{Eq:Coax_Faces} where $e$
 ranges over all unit vectors in $\R^3$.
 A nonempty intersection of two coaxial faces lies in the orbit of $SO(3)$.
 Each face $F_e^M$ is a circle of radius $\gamma$ and each face $F_e^m$ is a circle
 of radius $\alpha$.
 When $\chi$ has a repeated eigenvalue so that either $\alpha$ or $\gamma$ vanishes,
 then the corresponding coaxial face degenerates to a point.
\end{lemma}

A consequence of Lemma~\ref{L:Boundary_one} is that the coaxial faces are maximal faces,
and are therefore facets.\medskip

\noindent{\it Proof.}
 By Theorem~3.3 of~\cite{GLS05}, $V$ is the set of anisotropic tensors whose
 eigenvalues lie in the interval $[m,M]$, and so its boundary
 consists of tensors $\chibar$ either having maximal eigenvalue $M$
 or having minimal eigenvalue $m$.
 This shows that the boundary of $V$ consists of coaxial faces, which are thus the
 facets of $V$.

 We show that the intersection of two coaxial faces lies in the orbit of $\chi$.
 Suppose that $\chibar$ lies on two different coaxial faces.
 If these are $F^M_e$ and $F^m_f$, then $e$ and $f$ are eigenvectors of $\chi$
 with eigenvalues $M$ and $m$, respectively.
 The third eigenvalue of $\chibar$ is $-M-m$, and so $\chibar$ lies in the orbit of
 $\chi$.
 If the two faces have the form $F^M_e$ and $F^M_f$ with $e$ and $f$
 linearly independent, then the eigenvalue $M$ of $\chibar$ has multiplicity 2 and its
 third (smallest) eigenvalue is $-2M$.
 Since $\chibar\in V$, this smallest eigenvalue is bounded below by $m$; as
 $m\geq -2M$, we see that $m=-2M$ and so  again $\chibar$ lies in the orbit of
 $\chi$.
 The argument is similar if the two faces are $F^m_e$ and $F^m_f$.\smallskip

 The coaxial face $F^M_e$ of $V$ consists of tensors
 $\chibar\in V$ having $e$ as an eigenvector with eigenvalue $M$.
 Since each point of the boundary of $F_e^M$ lies in some other coaxial face,
 this boundary lies in the orbit of $\chi$ and is necessary an orbit of
 $Q_e$.
 We need only consider the case when $\chi\in F_e^M$ so that this boundary is
 $Q_e.\chi$.
 Suppose that $e=e_1$ and $\chi=\chi(M,0,0,\gamma,0)$, in the
 coordinates~\eqref{Eq:Coords}.
 Here, $\gamma=-\frac{M}{2}-m\geq 0$.
 As in Section~\ref{S:coaxial}, elements of $Q_e$ act on
 $\chi$ by rotation of the vector $(\gamma,0)$ formed by the last two coordinates,
 and thus $F^M_{e_1}$ is a circle of radius $\gamma$, which degenerates to a point
 if $\gamma=0$.

 We omit the similar arguments for $F^m_e$.
\QED\smallskip

\begin{theorem}
 If zero is not an eigenvalue of $\chi$, then $V$ has Carath\'eodory number $3$,
 and when zero is an eigenvalue, $V$ has Carath\'eodory number $2$.
\end{theorem}

When zero is not an eigenvalue of $\chi$, this is the main result about $V$
from~\cite{GLS05}.\medskip

\noindent{\it Proof.}
 First suppose that zero is not an eigenvalue of $\chi$.
 Every facet has Carath\'eodory number 2, as it is a circle.
 So by Lemma~\ref{L:Cartheodory}, $V$ has Carath\'eodory number either 2 or 3.
 By Lemma~\ref{L:Orbit_dimension}, 0 lies in $V$.
 If $V$ has Carath\'eodory number 2 then there exist $\lambda\in [0,1]$ and
 $R,S\in SO(3)$ with $0= \lambda R.\chi + (1-\lambda) S.\chi$.
 Multiplying by $R^{-1}$, this becomes $0 = \lambda\chi + (1-\lambda) R.\chi$, for a
 different rotation $R\in SO(3)$, and so
 \begin{equation}\label{Eq:CN=2}
    -\lambda\chi \ =\ (1-\lambda) R.\chi\,.
 \end{equation}
 Suppose that $\chi$ is diagonal.
 Then~\eqref{Eq:CN=2} implies that $R.\chi$ is also diagonal.

 If $M\neq -m$ so that 0 is not an eigenvalue of $\chi$, then one of the
 diagonal matrices $-\lambda\chi$ and  $(1-\lambda)R.\chi$ has two positive entries and
 the other has two negative entries, which is a contradiction.
 Thus if 0 is not an eigenvalue of $\chi$, then $V$ has Carath\'eodory
 number 3.\medskip

 Now we assume that 0 is an eigenvalue of $\chi$.
 We will show that the image of the map $[0,1]\times SO(3)\times SO(3)\to V$
 which takes $(\lambda,R,S)$ to $\lambda R.\chi+(1-\lambda)S.\chi$
 meets each $SO(3)$-orbit in $V$ and is therefore surjective.
 They key point is that two tensors are in the same orbit if and only if
 they have the same characteristic polynomial.

 The characteristic polynomial of a trace-zero matrix $\overline{\chi}$ with eigenvalues
 $s,t,-s-t$ is
\[
    x^3 - x(st+t^2+s^2) + (s^2t+st^2)\,.
\]
 The constant term is $-\det(\overline{\chi})$, while the coefficient $-\alpha$ of
 $x$ is the sum of the pairwise products of eigenvalues, which is an invariant of the
 matrix.

\begin{figure}[htb]
\[
  \begin{picture}(170,120)(-18,5)
   \put(0,11.2){\includegraphics[height=110pt]{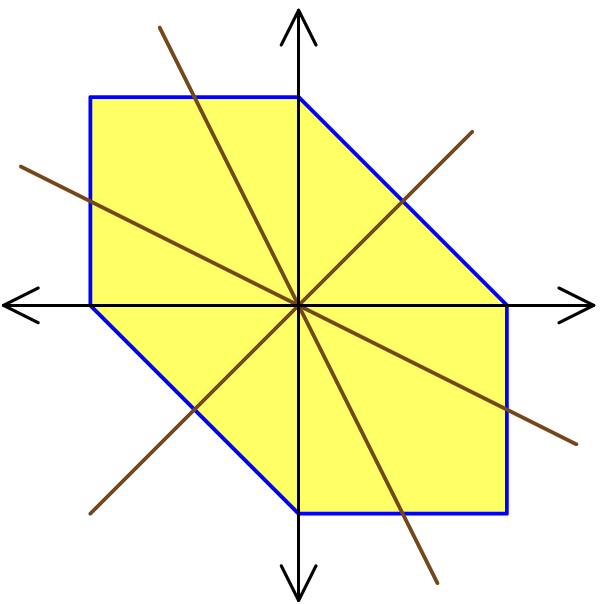}}
   \put(98,73){$1$}   \put(62,107){$1$}
   \put(-2,55){$-1$}  \put(36,17){$-1$}
   \put(87,108){\Brown{$t=s$}}
   \put(-19,119){\Brown{$t=-2s$}}
   \put(105,22){{\Brown{${\displaystyle t=-\frac{s}{2}}$}}}
  \end{picture}
    \qquad
  \begin{picture}(235,120)(0,-5)
   \put(0,0){\includegraphics{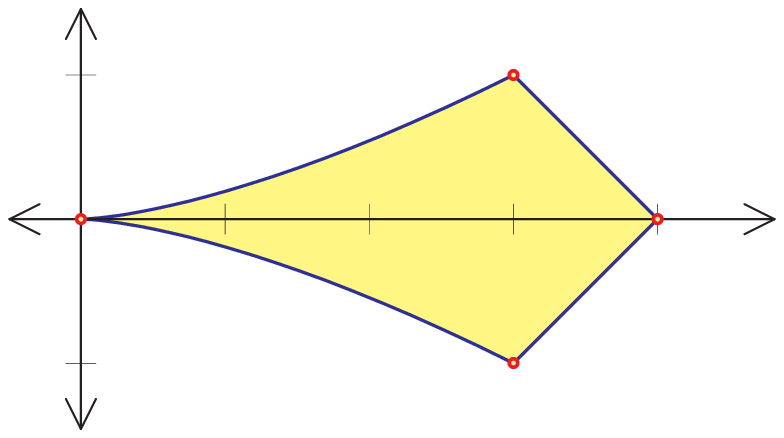}}

   \put(138,108){$(\frac{3}{4},\frac{1}{4})$}
   \put(133,  0){$(\frac{3}{4},-\frac{1}{4})$}
   \put(165,20){$\det=\alpha-1$}
   \put(165,90){$\det=1-\alpha$}
   \put(38, 0){$27\det^2=4\alpha^3$}
   \put(78,11){\vector(0,1){53}}
   \put(85,11){\vector(1,1){20}}
   \put(11,95){$\frac{1}{4}$}
   \put(2,12){$-\frac{1}{4}$}
   \put(189,42){$1$}
   \put(213,65){$\alpha$}
   \put(34,106){$\det$}
  \end{picture}
\]
\caption{Eigenvalues and invariants of tensors in $V$.}\label{Fig:Hex}
\end{figure}

 Scaling $\chi$, we may assume that its eigenvalues are $1,0$, and $-1$, so
 that $V$ consists of tensors $\overline{\chi}\in W$ with eigenvalues
 $s,t,-s-t$ lying in the interval $[-1,1]$.
 The set of such pairs $(s,t)$ are the points of the hexagon of Figure~\ref{Fig:Hex}.
The three lines through the origin $t=s$, $t=-s/2$, and $t=2s$ divide the hexagon into six
quadrilaterals and permutations of the eigenvalues permute these
quadrilaterals.
We leave it to the reader to check that $(s,t)\mapsto(\alpha,\det)$ is a one-to-one
mapping of each quadrilateral onto the region
shown in Figure~\ref{Fig:Hex}, which is
 \begin{equation}\label{Eq:Region}
   \Blue{X}\ :=
   \left\{(\alpha,\det)\in\R^2\mid 27{\textstyle \det^2}\leq 4\alpha^3 \mbox{\ and\ }
               \alpha\leq 1-|\det|\right\}\,,
 \end{equation}
 and is bounded by the curves  $\det=1-\alpha$,
 $\det=\alpha-1$, and $27\det^2=4\alpha^3$.

Consider matrices of the form
$\chi(\lambda,\theta,\tau):=\lambda R(\theta).\chi + (1-\lambda)S(\tau).\chi$, where
\[
  \Blue{R(\theta)}\ :=\ \left[
   \begin{matrix}\cos\theta&\sin\theta&0\\
                 -\sin\theta&\cos\theta&0\\
                 0&0&1\end{matrix}\right]
  \qquad\mbox{and}\qquad
  \Blue{S(\tau)}\ :=\ \left[
   \begin{matrix}\cos\tau&0&\sin\tau\\0&1&0\\
                -\sin\tau&0&\cos\tau\end{matrix}\right]\,.
\]
%
%
The invariants $(\alpha,\det)$ of  $\chi(\lambda,\theta,\tau)$ are
\[
   \left(1- \lambda(1-\lambda)(4\sin^2\tau + \sin^2\theta - 2\sin^2\tau\sin^2\theta),
   \quad
   \lambda(1-\lambda)\sin^2\theta(1-2\lambda\sin^2\tau)\right)\,.
\]
%
If we let $u=\sin^2\theta$ and $v=\sin^2\tau$, then the set of invariants
of  $\chi(\lambda,\theta,\tau)$ for all $(\lambda,\theta,\tau)$ are the image of the unit
cube $[0,1]^3$ under the map
\[
  f\ \colon\ (\lambda,u,v)\ \longmapsto\
   \left( 1-\lambda(1-\lambda)(4v+u-2uv), \ \lambda(1-\lambda)u(1-2\lambda v) \right)\,.
\]

We show that the image of $f$ includes that part of $X$~\eqref{Eq:Region} where
$\det\geq 0$.
This will complete the proof, as replacing $\chi$ by $-\chi=S(\frac{\pi}{2}).\chi\in V$
in our definition of $f$ changes the sign of the determinant and does not change the
invariant $\alpha$.

\begin{figure}[htb]
\[
  \begin{picture}(180,105)(0,0)
  \put(0,0){\includegraphics[height=90pt] {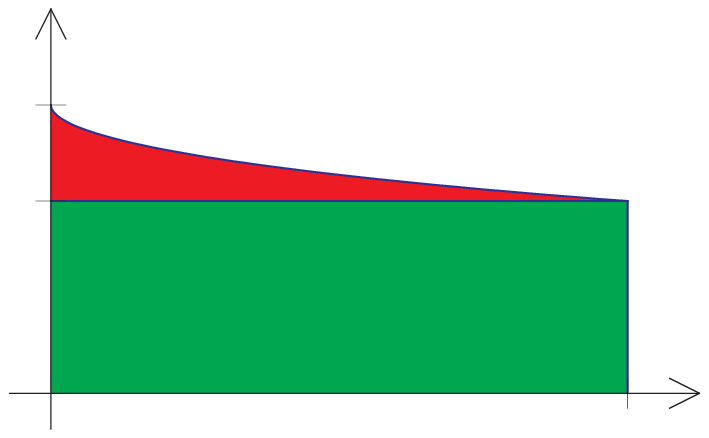}}
    \put(44,80){$(\lambda-\lambda^2)u=  2(1-2\lambda)^2$}
   \put(23,100){$\lambda$}  \put(168,20){$u$}
   \put(2,57){$\frac{1}{3}$}   \put(2,80.8){$\frac{1}{2}$}
   \put(156,-3){$1$}
   \put(60,-10){Face $v=1$}
  \end{picture}
  \qquad\qquad
  \begin{picture}(180,105)(0,0)
  \put(0,0){\includegraphics[height=90pt] {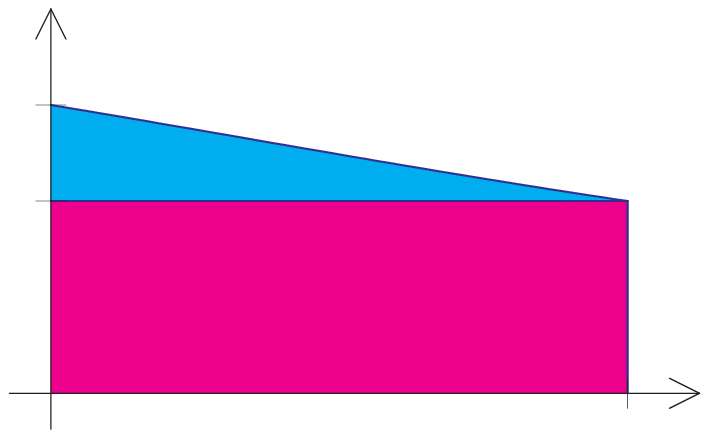}}
   \put(44,84){$2(\lambda-\lambda^2)v = -(1-\lambda-2\lambda^2)$}

   \put(23,100){$\lambda$}  \put(168,20){$v$}
   \put(2,57){$\frac{1}{3}$}   \put(2,80.8){$\frac{1}{2}$}
   \put(156,-3){$1$}
   \put(60,-10){Face $u=1$}
  \end{picture}
\]
\caption{Subsets of faces of the cube}\label{F:subsets_faces}
\end{figure}

\begin{figure}[htb]
\[
  \begin{picture}(335,125)(-23,-32)
   \put(0,0){\includegraphics[height=70pt]{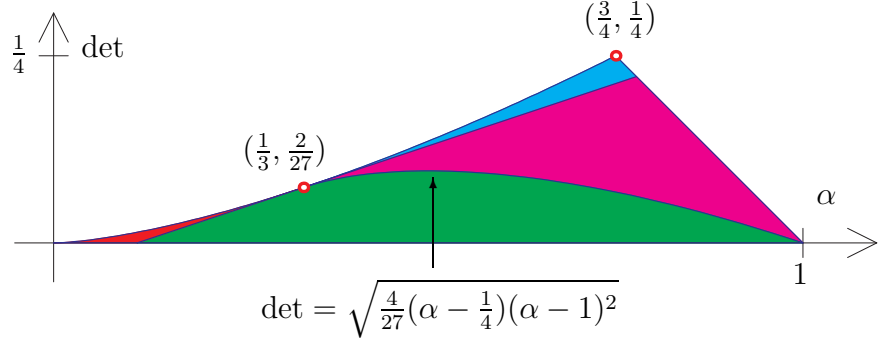}}
   \put(-15,69){$\frac{1}{4}$}   \put(281,-15){$1$}
   \put( 73,32){$(\frac{1}{3},\frac{2}{27}$)}
   \put(202,82){$(\frac{3}{4},\frac{1}{4}$)}
   \put(80,-28){$\det=\sqrt{\frac{4}{27}(\alpha-\frac{1}{4})(\alpha-1)^2}$}
   \put(145,-9){\vector(0,1){34}}
   \put(12,71){$\det$} \put(290,15){$\alpha$}
  \end{picture}
\]
\caption{Image of subsets of Figure~\ref{F:subsets_faces}}\label{F:images}
\end{figure}

Figure~\ref{F:subsets_faces} shows subsets of the faces $v=1$ and
$u=1$ of the cube
 which include the segments $\lambda=1/3$ and have boundaries the indicated curves.
 The map $f$ is one-to-one on the interior of each region, and the
 images cover that part of $X$ with $\det\geq 0$, meeting only along the curve
 $\det=\sqrt{\frac{4}{27}(\alpha-\frac{1}{4})(\alpha-1)^2}$
 for $\frac{1}{3}\leq \alpha\leq 1$, as shown in Figure~\ref{F:images}.
 The line $3\det=(\alpha-\frac{1}{9})$, which is tangent to the
 boundary curve $27\det^2=4\alpha^3$ at the point $(\frac{1}{3},\frac{2}{27})$,
 is the image of the lines $\lambda=\frac{1}{3}$ in Figure~\ref{F:subsets_faces}.
 This completes the proof.
\QED\smallskip

 Lastly, we identify the hyperplanes supporting facets of $V$.
 By Lemma~\ref{L:Boundary_one}, the faces $F_e^M$ and $F_e^m$ are 2-dimensional, unless
 $\chi$ has an eigenvalue with multiplicity 2, and in that case exactly one face is
 degenerate.

 \begin{proposition}
  The facets $F_e^M$ and $F_e^m$ have a unique supporting linear function, unless they are
  degenerate.
 \end{proposition}

\noindent{\it Proof.}
 It suffices to determine the hyperplanes supporting faces which contain $\chi$.
 Any hyperplane supporting the vertex $\chi$ contains the
 the tangent space $T_\chi\calO$ at
 $\chi$ to the orbit $\calO$ through $\chi$.
 Choose coordinates so that
 $\chi=\chi(M,0,0,\gamma,0)$ is diagonal ($\gamma=-\tfrac{M}{2}-m$) so that
 $\chi\in F^M_{e_1}\cap F^m_{e_3}$.
 Recall that $T_\chi\calO$~\eqref{Eq:TanSpace_k=1} is the affine 3-plane in $W$
\[
   \chi(M, \Magenta{*}, \Magenta{*}, \gamma, \Magenta{*})\,,
\]
 where $\Magenta{*}$ represents an arbitrary real number.
 If $(x_1,\dotsc,x_5)$ are the coordinates~\eqref{Eq:Coords}, then hyperplanes
 containing $T_\chi\calO$ have equation $\sum_i c_ix_i=c$, where
 \begin{equation}\label{Eq:Supp_1}
  c_2\ =\ c_3\ =\ c_5\ =\ 0,\qquad\mbox{and}\qquad
   c_1M + c_4 \gamma\ =\ c\,.
 \end{equation}
 If $\gamma\neq 0$, then $F^M_{e_1}$ is nondegenerate and contains the additional
 point
\[
   R_{e_1,\frac{\pi}{2}}\ =\ \chi(M,0,0,0,\gamma)\ =\
   \left(\begin{matrix}M&0&0\\0&-\tfrac{M}{2}&\gamma\\
                       0&\gamma&-\tfrac{M}{2}\end{matrix}\right)\ ,
\]
 which imposes the further condition $c_1 M = c$ on a support hyperplane to
 $F^M_{e_1}$.
 Thus $c_1=c/M$ and $c_4=0$.
 Setting $c=M$ so that $c_1=1$, we see that the support hyperplane to
 $F^M_{e_1}$ is defined by $x_1 = M$, which is $L_{e_1}(\chi)=M$.

 If $\nu=M+\frac{m}{2}\neq 0$, then $F^m_{e_3}$ is nondegenerate and it
 contains the point
\[
   \left(\begin{matrix}-\tfrac{m}{2}&\nu&0\\\nu&-\tfrac{m}{2}&0\\0&0&m\end{matrix}\right)
   \ =\ \chi(-\tfrac{m}{2},\nu,0,-\tfrac{3m}{4},0)\,,
\]
 and so a support hyperplane to $F^m_{e_3}$ must satisfy~\eqref{Eq:Supp_1} and also
 $-c_1\tfrac{m}{2} - c_4 \tfrac{3m}{4} = d$.
 Subtracting these equations and dividing by $\nu$, we see that $c_1=2c_4$, and
 so
\[
  x_1 + 2x_4\ =\ -2m\,,
\]
 is the support hyperplane to $F^m_{e_3}$.
 Note that $x_1+2x_4$ is $-2L_{e_3}(\chi)$.
\QED

%
%
\section{Two Metals}\label{S:Two}
%

Let $\chi_1,\chi_2\in W$ be linearly independent anisotropic tensors
and set $\chi:=(\chi_1,\chi_2)\in W^2$.
By Lemma~\ref{L:Orbit_dimension}, the convex hull $V$ of the orbit
$\calO:=SO(3).\chi$ is a $2\cdot 5=10$-dimensional convex body containing the origin.
Its boundary is 9-dimensional.

We study the facial structure of $V\subset W^2$.
One tool will be a family of $SO(3)$-equivariant maps $\pi_\alpha\colon W^2\to W$.
We first determine the dimension of the orbit, show that
the maximum dimension of a facet is 6, and then define coaxial faces.
Our main result is that coaxial faces are facets if $\chi_1$ and
$\chi_2$ have distinct eigenvectors.
In that case, almost all coaxial faces have dimension 6 and Carath\'eodory number
4.
We are unable to rule out the existence of other facets, but we conjecture that there are
no other facets.

Let $\Span(\chi)\subset W$ be the 2-dimensional subspace of $W$ spanned by $\chi_1$ and
$\chi_2$.
The structure of $V$ depends only on $\Span(\chi)$.
Indeed, if $\chi'_1,\chi'_2\in\Span(\chi)$ are linearly independent, then
there is a $2\times 2$ invertible matrix $A=(a_{ij})$ such that
\[
   \chi'_1\ =\ a_{11}\chi_1\;+\; a_{12}\chi_2
       \qquad
   \chi'_2\ =\ a_{21}\chi_1\;+\; a_{22}\chi_2\,.
\]
This induces an $SO(3)$-isomorphism $W^2\xrightarrow{\sim}W^2$:
\[
   W^2\ \ni\ (w_1,w_2)\ \longmapsto\ ( a_{11}w_1\;+\; a_{12}w_2,\
            a_{21}w_1\;+\; a_{22}w_2)\ \in\ W^2
\]
which sends $V$ to the convex hull of the orbit of
$(\chi'_1,\chi'_2)$.
This is nothing more than a change of coordinates on $W^2$.

Any non-zero vector $\alpha=(\alpha_1,\alpha_2)\in\R^2$ gives an
$SO(3)$-map
 \begin{equation}\label{Eq:pi_alpha}
   \pi_\alpha\ \colon\ W^2\ \longrightarrow\ W
 \end{equation}
defined by $\pi_\alpha(w_1,w_2):=\alpha_1w_1+\alpha_2w_2$.
Write $w_\alpha$ for $\pi_\alpha(w)$.
In particular, $\Blue{\chi_\alpha}:=\pi_\alpha(\chi)\in\Span(\chi)$.
Set $\Blue{V_\alpha}\subset W$ to be the convex hull of the orbit
$\Blue{\calO_\alpha}:=SO(3).\chi_\alpha$.
Since $\pi_\alpha(SO(3).\chi)=SO(3).\chi_\alpha$, we have
$\pi_\alpha(\calO_\chi)=\calO_\alpha$ and $V_\alpha=\pi_\alpha(V)$.
We compute the dimension of the orbit $\calO_\chi$.

\begin{theorem}
 $\dim \calO_\chi =3$.
\end{theorem}

\noindent{\it Proof.}
 We will show that $\dim\calO_\alpha=3$ for some $\alpha\in\R^2$.
 As $\pi_\alpha(\calO_\chi)=\calO_\alpha$, this implies $\dim\calO_\chi\geq 3$.
 Since $\dim SO(3)=3$, we have $\dim\calO_\chi\leq 3$
 and so $\dim \calO_\chi =3$.

 By Proposition~\ref{Prop:orbit_dim},  the dimension of $\calO_\alpha$ is
 3 if and only if $\chi_\alpha$ has distinct eigenvalues.
 If either $\chi_1$ or $\chi_2$, say $\chi_1$, has distinct eigenvalues,
 then $\dim\calO_{(1,0)}=3$
 and we are done.
 Suppose the contrary, that neither $\chi_1$ nor $\chi_2$ has distinct eigenvalues.
 That is, for each $i=1,2$, $\chi_i$ has a 2-dimensional eigenspace with eigenvalue
 $\alpha_i$.
 Since $0$ cannot be a repeated eigenvalue, neither $\alpha_1$ nor
 $\alpha_2$ is zero.
 These eigenspaces must meet, so $\chi_1$ and $\chi_2$ share an eigenvector,
 which is an eigenvector for the nonzero tensor
 $\chi_\alpha:=\alpha_2\chi_1\Red{-}\alpha_1\chi_2\in\Span(\chi)$
 with eigenvalue 0.
 But then $\chi_\alpha$ has distinct eigenvalues and so $\dim\calO_\alpha=3$.
\QED\smallskip

\begin{lemma}
 The maximum dimension of a proper face of $V$ is $6$.
\end{lemma}

\noindent{\it Proof.}
Let $F$ be a proper face of $V$ and let $S\subset SO(3)$ its stabilizer subgroup,
\[
  S\ =\  \left\{g\in SO(3)\mid g.F\subset F\right\}\,.
\]
This is a closed, proper subgroup, and thus either has dimension 1
(in which case it is a coaxial subgroup $Q_e$ or a rotation subgroup $Q_e^+$), or
it is finite and has dimension zero.

Let $F^\circ$ be the relative interior of $F$, those points of
$F$ which do not lie in any other face of $V$ of the same or smaller
dimension. If $g.F^\circ\cap F^\circ\neq \emptyset$ then $g.F=F$,
and so $g\in S$.

Let \Blue{$\partial V$} be the boundary of $V$
and consider the map $f\colon SO(3)\times F^\circ\to \partial V$ defined by
 \[
   f\ \colon\ (g,v)\ \longmapsto\  g.v\,.
 \]
This map is not 1-1:  Suppose that $g.v=h.w$, for $g,h\in SO(3)$ and
$v,w\in F^\circ$. Then $h^{-1}g.v=w$ and so  $h^{-1}g.F^\circ\cap
F^\circ\neq \emptyset$, which implies that $s:=h^{-1}g\in S$. Then
$s.v=w$.

This calculation shows that the fibers of $f$ have the form
\[
   \left\{ (gs^{-1}, sv)\mid s\in S\right\}, \qquad\mbox{for}\ g\in SO(3)
    \ \mbox{and}\ x\in F^\circ\,.
\]
Thus we have the dimension calculation
\[
   \dim \partial V\ \geq\ \dim SO(3)+\dim F^\circ-\dim S\,.
  \]
Since $\dim \partial V=9$ and $\dim SO(3)=3$, this gives
 \begin{equation}\label{eq:dimension_ineq}
  6+\dim S\ \geq\ \dim F\,.
 \end{equation}
If $S$ is finite, then $\dim F\leq 6$.
If $S$ has dimension 1 so that it is either $Q_e$ or $Q_e^+$ for some $e$,
then $F$ could have dimension up to 7.
By~\eqref{Eq:Decomposition}, $W^2=\R^2\oplus U_1^2\oplus U_2^2$ as a
representation of $S$.
If $S=Q_e$, then Lemma~\ref{L:Orbit_dimension} implies that $F$ has even dimension,
and if $S=Q_e^+$, then Example~\ref{Ex:SO(2)_Orbits} implies that $F$ has dimension 0, 2,
or 4, which completes the proof.
\QED

%
\subsection{Coaxial faces}
%

A \Blue{{\it coaxial face}} of $V$ is a face that is stabilized by some
coaxial subgroup, $Q_e$.
By Lemma~\ref{Lem:trivial_support}, such a face is supported by a $Q_e$-invariant linear
function, which must factor through the projection to the
trivial isotypical component of $W^2$, by Schur's Lemma.
Since this component is $\R^2$~\eqref{Eq:Decomposition}, $L$ is the pullback of
a linear map 
\[
   \R^2\ni(M_1,M_2)\ \longmapsto\ \alpha_1 M_1+\alpha_2 M_2\in \R\,.
\]
Up to a scalar, this is the composition of the $Q_e$-invariant linear function
$L_e$~\eqref{E:L_e} on $W$ with $\pi_\alpha$, which is the map
$\Blue{L_{e,\alpha}}$ defined by
\[
  L_{e,\alpha}(w)\ :=\ L_e(w_\alpha)\ =\ \langle e, w_\alpha e\rangle\,.
\]

Suppose now that $e$ is a unit vector.
For each non-zero $\alpha\in\R^2$, define
 \begin{eqnarray*}
  M_\alpha&:=& \mbox{maximum eigenvalue of }\chi_\alpha\,,\quad\mbox{and}\\
  m_\alpha&:=& \mbox{minimum eigenvalue of }\chi_\alpha\,.
 \end{eqnarray*}
 If $\chibar\in V$, then $\chibar_\alpha\in V_\alpha$, and
 so by Lemma~\ref{Lem:Le-bounds} we have
 \begin{equation}\label{Eq:Lambda_alpha}
   M_\alpha\ \geq\ L_{e,\alpha}(\chibar)\ \geq\ m_\alpha\,,
 \end{equation}
 with equality only when $e$ is an eigenvector of $\chibar_\alpha$ having
 eigenvalue $M_\alpha$ or $m_\alpha$.
 Thus coaxial faces are the faces of $V$ defined by equality in~\eqref{Eq:Lambda_alpha}.

For $e\in\R^3$ a unit vector and $0\neq\alpha\in\R^2$, define the coaxial face
\[
   \Blue{F_{e,\alpha}}\ :=\ \left\{\chibar\in V\mid L_{e,\alpha}(\chibar)=M_\alpha\right\}\,.
\]
 If $\chibar\in F_{e,\alpha}$, then $e$ is an eigenvector of $\chibar_\alpha$ with
 eigenvector $M_\alpha$.
 As in Section~\ref{S:One}, each coaxial face $F_{e,\alpha}$ is the convex hull of an
 orbit $Q_e.\chi'$, for some $\chi'\in\calO$.

\begin{theorem}\label{T:coaxial}
  Faces of $V$ have dimension at most $6$.
  The coaxial faces of\/ $V$ form a $3$-dimensional family
  whose union is a $9$-dimensional
  subset of the boundary of\/ $V$ if and only if $\chi_1$ and $\chi_2$ have distinct
  eigenvectors.
  When this happens, almost all coaxial faces have dimension $6$, have Carath\'eodory
  number $4$, and are facets of\/ $V$.
\end{theorem}

If the boundary of $V$ is the union of the coaxial faces, then Lemma~\ref{L:Cartheodory}
implies that the Carath\'eodory number of $V$ is at most $5$, and we conjecture this is
the case.
If there are faces of dimension 6 that are not coaxial, then Carath\'eodory's Theorem
implies that their Carath\'eodory number is at most 7.
Then Lemma~\ref{L:Cartheodory} implies the following corollary of Theorem~\ref{T:coaxial}.

\begin{cor}\label{co:Car}
 The Carath\'eodory number of $V$ is at most $8$.
\end{cor}

By almost all in the statement of Theorem~\ref{T:coaxial}, we mean in the algebraic sense:
Except for those $\alpha\in\R^2$ lying in finitely many half-rays in $\R^2$,
$F_{e,\alpha}$ has dimension 6 when $\chi_1$ and $\chi_2$ have distinct eigenvectors.
The proof of Theorem~\ref{T:coaxial} is done in the series of lemmas below.

\begin{remark}
  The condition that the magnetic susceptibility tensors $\chi_1$ and $\chi_2$ have
  distinct eigenvectors has already been considered in protein folding.
  It implies that RDC measurements from the two ions are sufficient to
  remove the symmetry property of the RDC~\cite{LPS}.
\end{remark}

 Since $M_{-\alpha}=-m_\alpha$, there is no need for two types of coaxial faces as in
 Section~\ref{S:One}.
 Since if $r>0$ then $M_{r\alpha}=rM_\alpha$ and $L_e=L_{-e}$, we have
\[
   F_{e,\alpha}\ =\ F_{-e,\alpha}\ =\ F_{e,r\alpha}
\]
 if $r>0$.
 Thus we may assume that $\alpha$ lies on the unit circle $S^1$ in $\R^2$.
 We also only need to consider the unit vector $e$ up to multiplication by $\pm 1$,
 that is, as a point in the real projective plane, $\R\P^2$,
 which is a 2-dimensional manifold.

\begin{lemma}\label{Lem:Coaxial_equivariant}
 The coaxial faces $F_{e,\alpha}$ form a $3$-dimensional family parameterized
 by $\R\P^2\times S^1$.
 For each $\alpha\in S^1$, any two coaxial faces $F_{e,\alpha}$ and $F_{e',\alpha}$ are
 isomorphic.
\end{lemma}

Since the boundary of $V$ is 9-dimensional and it has a 3-dimensional family of coaxial faces,
we see again that the maximum dimension of a coaxial face is 6.\smallskip

\noindent{\it Proof.}
 Suppose that $e\in\R^3$ is an eigenvector for $\chi_\alpha$ with maximal
 eigenvalue $M_\alpha$.
 Then $\chi\in F_{e,\alpha}$ and $F_{e,\alpha}$ is the convex hull of the orbit $Q_e.\chi$.
 If $R\in SO(3)$, then
 \[
    R.Q_e.\chi\ =\  RQ_eR^T. R.\chi\ =\  Q_{Re}. (R.\chi)\,.
 \]
 But $(R.\chi)_\alpha$ is an anisotropic tensor having eigenvector $Re$ with eigenvalue
 $M_\alpha$.
 Therefore $F_{Re,\alpha}$ is the convex hull of $Q_{Re}. (R.\chi)=R.Q_e.\chi$, and
 thus equals $R.F_{e,\alpha}$.
\QED\smallskip

We now determine the dimension of the coaxial faces. By
Lemma~\ref{Lem:Coaxial_equivariant}, we need only study one coaxial
face $F_{e,\alpha}$ for each $\alpha\in S^1$. We compute the
dimension of the affine span of an orbit $Q_e.\chi$, where $e$ is an
eigenvector of $\chi_\alpha$. This is the dimension of a coaxial
face when the eigenvalue associated to $e$ is a maximal eigenvalue
of $\chi_\alpha$. Since this dimension is the rank of a matrix,
those entries are algebraic functions of $\alpha$. Thus that for all
but finitely many $\alpha$, this rank will be constant and it will
be smaller for $\alpha$ in that finite set.

Let $\alpha\in S^1$ and suppose that $\chi_\alpha,\chi'\in\Span(\chi)$ are linearly
independent, and let $e$ be a unit eigenvector of $\chi_\alpha$.
The decomposition~\eqref{Eq:Decomposition} of $W$ into $Q_e$-isotypical components induces
a decompostion of the tensors $\chi_\alpha$ and $\chi'$ into
their components in $\R\oplus U_1\oplus U_2$,
 \[
   \chi_\alpha\ =\ M_\alpha \oplus 0  \oplus y_\alpha
    \qquad\mbox{and}\qquad
   \chi'      \ =\ M'       \oplus x' \oplus y'      \,.
 \]
The $U_1$-component of $\chi_\alpha$ is $0$, because $e$ is an
eigenvalue of $\chi_\alpha$. Let $d_1\in\{0,1\}$ be the dimension of
the linear span of $x'$ in $U_1$ and $d_2\in\{0,1,2\}$ be the
dimension of the linear span of $y_\alpha,y'$ in $U_2$. By
Lemma~\ref{L:Orbit_dimension}, the dimension of the convex hull of
$Q_e.\chi$ is $d_1\cdot\dim U_1 + d_2\cdot\dim U_2$, which implies
the following lemma.

\begin{lemma}\label{L:face_dim}
  The coaxial face $F_{e,\alpha}$ has dimension $2(d_1+d_2)$.
\end{lemma}

Thus again a coaxial face has dimension at most 6.

\begin{lemma}
 If $\chi_1$ and $\chi_2$ have a common eigenvector, then coaxial faces have dimension
 $2$ or $4$.
\end{lemma}

\noindent{\it Proof.}
 Fix $\alpha\in S^1$ and let $f$ be a common eigenvector of $\chi_1$ and $\chi_2$.
 Then it is an eigenvector of any $\chi_\alpha$.
 Let $\chi'$ be another tensor in $\Span(\chi)$ which is not proportional to
 $\chi_\alpha$.
 Suppose that $e=e_1$ is an eigenvector of $\chi_\alpha$, that $f\in\{e_1,e_2,e_3\}$, and
 write $\chi_\alpha$ and $\chi'$ in the coordinates~\eqref{Eq:Coords},
\[
    \chi_\alpha\ =\ (M_\alpha,\ 0,0,\ y,z)
   \qquad\mbox{and}\qquad
   \chi'\ =\ (M',\ w',x',\ y',z')\,.
\]
 Note that $(w',x')\in U_1$ and $(y,z), (y',z')\in U_2$.

 If $f=e_1$, then $(w',x')=(0,0)$, and so $d_1=0$.
 If $f=e_2$ or $e_3$, then $z=z'=0$ and so $d_2=1$.
 In either case, $d_1+d_2<3$ and so the coaxial face $F_{e,\alpha}$ has dimension $2$ or
 $4$.
\QED

\begin{lemma}
  If $\chi_1$ and $\chi_2$ do not have a common eigenvector, then there is a
  coaxial face with dimension $6$.
\end{lemma}

\noindent{\it Proof.}
  Suppose that $e=e_1$ and $e_1,e_2,e_3$ is an ordered basis of eigenvectors of
  $\chi_1$ with the eigenvalue of $e$ maximal.
  Write $\chi_1$ and $\chi_2$ in the coordinates~\eqref{Eq:Coords},
\[
  \chi_1\ =\ \chi(M_1,\ 0,0,\  \gamma,0)
   \qquad
  \chi_2\ =\ \chi(M_2,\ w,x,\ y,z)\,.
\]

 The dimension of the coaxial face $F_{e_1,(1,0)}$ is 6 if and only if
 $(w,x)\neq (0,0)$ and $(\gamma,0), (y,z)\in U_2$ are linearly independent.
 Suppose that $\dim F<6$.
 We cannot have $(w,x)= (0,0)$ for then $e_1$ is a common eigenvector of $\chi_1$ and
 $\chi_2$, a contradiction.
 Thus the vectors $(\gamma,0), (y,z)$ are dependent.

 If $\gamma=0$, then $\chi_1$ has a repeated smallest eigenvalue with eigenspace
 spanned by $e_2$ and $e_3$.
 Changing the last two coordinates, we may assume that $z=0$.
 We cannot also have either $w=0$ or $x=0$ for then $\chi_1$ and $\chi_2$ have
 either $e_2$ or $e_3$ as a common eigenvector.
 If $y=0$, then $-xe_2+we_3$ is a common eigenvector, so $y\neq 0$ and $F_{e_1,(1,0)}$ has
 dimension 4.
 In the coordinates~\eqref{Eq:Coords} with respect to the ordered basis $e_3,e_2,e_1$,
 $-\chi_1$ and $\chi_2$ are
\[
  -\chi_1\ =\ \chi(\tfrac{M_1}{2},\, 0,0,\, 3\tfrac{M_1}{4},0)
   \quad{\rm and }\quad
  \chi_2\ =\
  \chi(-\tfrac{M_2}{2}-y,\,0,x,\,-3\tfrac{M_2}{4}+\tfrac{y}{2},w)\,.
\]
 Since $w,x\neq 0$, the affine span of $Q_{e_3}.\chi$ has dimension 6.
 Since $\tfrac{M_1}{2}$ is the maximal eigenvalue of $-\chi_1$ with
 eigenvector $e_3$, this shows that $F_{e_3,(-1,0)}$ has dimension 6.

 The third possibility is that $z=0$.
 But then the same arguments as in the previous paragraph
 show that coaxial face $F_{e_3,(-1,0)}$ has dimension 6.
\QED

%
\subsection{Structure and Carath\'eodory number of a coaxial facet}
%

Suppose that $F$ is a coaxial face of dimension 6.
We may assume that $F$ is the convex hull of the orbit $Q_e.\chi$ and that $F$ spans the
representation $U_1\oplus U_2^2$.
This $Q_e$-orbit is the union of two orbits of its identity component
$Q_e^+ (\simeq SO(2))$.
Call them $\calO^+$ and $\calO^-$.
By Example~\ref{Ex:SO(2)_Orbits}, each orbit spans a subrepresentation of
$U_1\oplus U_2^2$ isomorphic to $U_1\oplus U_2$.
Set
\[
   \Blue{F^{\pm}}\ :=\ \mbox{convex hull of }\calO^{\pm}
    \quad\mbox{ and }\quad
   \Blue{W^{\pm}}\ :=\ \mbox{linear span of }\calO^{\pm}\simeq U_1\oplus U_2\,.
\]

\begin{proposition}\label{Pr:Fpm}
 The faces $F^{\pm}$ each have dimension $4$ and Carath\'eodory number $3$.
 Points on their boundary are the convex hull of one or two vertices, while
 points in their relative interiors are the convex hull of three vertices.
\end{proposition}

\noindent{\it Proof.}
 As in Example~\ref{Ex:SO(2)_Orbits}, we identify $Q_e^+\simeq SO(2)$ with the circle
 group $S^1$ and $U_1,U_2$ with $\C$.
 Then $z\in S^1$ acts on $U_1$ as scalar multiplication by $z$ and on $U_2$ as scalar
 multiplication by $z^2$ and $F^{\pm}$ has dimension 4.

 Let $(u,v)^T\in \C^2\simeq U_1\oplus U_2$ be the point corresponding to $\chi$.
 Then
\[
  \calO^+\ =\ \left\{ (e^{i\theta} u, e^{2i\theta} v)\mid 0\leq \theta<2\pi\right\}\,,
\]
 and its convex hull is
\[
  \Bigl\{ \Bigl( \sum_{j=1}^n \lambda_j e^{i\theta_j} u,
       \sum_{j=1}^n \lambda_j e^{2i\theta_j} v\Bigr)\mid \sum_j \lambda_j=1,\
       \ 0\leq \theta_1,\dotsc,\theta_n< 2\pi\Bigr\}\,.
\]
But this is $B.(u,v)^T$, where $B$ is the set of $2\times 2$ diagonal matrices
whose entries are
 \begin{equation}\label{Eq:CF}
   \Bigl( \sum_{j=1}^n \lambda_j e^{i\theta_j} ,
       \sum_{j=1}^n \lambda_j e^{2i\theta_j}\Bigr)
   \quad\mbox{ where }\quad \sum_j \lambda_j=1
   \quad\mbox{ and }\quad
   0\leq \theta_1,\dotsc,\theta_n< 2\pi\,.
 \end{equation}
Thus $F^+$ (and also $F^-$) is isomorphic to $B$.

Curto and Fialkow~\cite{CF} characterized the points of $B$.
Let $(a,b)$ be a point of $B$~\eqref{Eq:CF}
and $p$ the corresponding measure on $S^1$,
\[
   p(e^{i\theta})\ =\ \left\{
    \begin{array}{lcl}  \lambda_j&\ &\mbox{if }\theta=\theta_j\,,\\
         0&& \mbox{otherwise}\,.\end{array}\right.
\]
Set $\gamma_{ij}:=\int_{S^1}\overline{z}^i z^jdp(z)$ for
$0\leq i,j$ with $i+j\leq 2$ and form the moment matrix
\[
   M\ :=\ \left(\begin{matrix}
          \gamma_{00}&\gamma_{01}&\gamma_{10}\\
          \gamma_{10}&\gamma_{11}&\gamma_{20}\\
          \gamma_{01}&\gamma_{02}&\gamma_{11}\end{matrix}\right)
   \quad =\quad
          \left(\begin{array}{ccc} 1&a&\overline{a}\\
             \overline{a}&1&\overline{b}\\a&b&1\end{array}\right)\ .
\]

\begin{proposition}[Curto and Fialkow~\cite{CF}]\label{Prop:CF}
  The points $(a,b)\in B$ are exactly the points $(a,b)\in\C^2$ such that $M$ is positive
  semi-definite.
  The rank of $M$ is the minimum number of summands needed to represent the point
  $(a,b)$~$\eqref{Eq:CF}$.
\end{proposition}

In particular, this implies that each body $F^\pm$ has Carath\'eodory number 3.

By Proposition~\ref{Prop:CF},
\[
  B\ =\ \left\{(a,b)\in\C^2\mid  1-|a|^2\geq 0,\ 1-|b|^2\geq 0,\
     1+\overline{a}^2b+a^2\overline{b}-2|a|^2-|b|^2\geq 0\right\}\,.
\]
If $1-|a|^2=0$ then $a\in S^1$ and so $n=1$ in~\eqref{Eq:CF}. If
$1-|b|^2=0$, then $b\in S^1$ and either $n=1$ in~\eqref{Eq:CF} or
$n=2$ with $|\theta_1-\theta_2|=\pi$. Thus if $M$ has rank 3, then
$|a|<1$ and $|b|<1$, and so $(a,b)$ lies in the interior of $B$ as
the inequalities are strict. This implies that points on the
boundary of $B$ are the convex hull of one or two vertices and this
completes the proof of Proposition~\ref{Pr:Fpm}. \QED\smallskip

Now we complete the proof of Theorem~\ref{T:coaxial}, showing that
the coaxial facet $F$ has Carath\'eodory number 4. The coaxial face
$F$ is the convex hull of $F^+$ and $F^-$. Let $v\in F$. We suppose
that $v\not\in F^+\cup F^-$, for otherwise $v$ is the convex
combination of at most three vertices. Then there exist $v^\pm\in
F^\pm$ and $\lambda\in(0,1)$ such that
\[
   v\ =\ \lambda v^+ + (1-\lambda)v^-\,.
\]
If both $v^+$ and $v^-$ lie on the boundary of their respective
subfaces, then each is a convex combination of at most 2 vertices,
and $v$ is a convex combination of at most 4 vertices.

Suppose instead that $v^+$ lies in the relative interior of $F^+$.
The linear span of $v$ and $W^-$ has dimension 5 in the 6-dimensional space
$U_1\oplus U_2^2$ and therefore it meets $W^+$ in a 3-dimensional affine subspace $U^+$.
Similarly the span of $v$ and $W^+$ meets $W^-$ in a 3-dimensional subspace $U^-$.
Observe that both $U^+$ and $U^-$ contain the two-dimensional linear subspace
$W^+\cap W^-=U_1$ so that their span has dimension 4.

Consider the cone over $U^-\cap F^-$ with vertex $v$.
Removing $v$, this has two components.
One meets $F^-$.
Let $C$ be the  component which does not meet $F^-$, and let
$\Blue{C^+}:=C\cap U^+$.
This is a convex set which contains $v^+$ and thus meets the relative interior of
$\Blue{B^+}:=U^+\cap F^+$.
Points $v'\in C^+\cap B^+$ are exactly those points of $F^+$ for which there exists
a point $v''\in F^-$ such that $v$ is a convex combination of $v'$ and $v''$.
There are two possibilities.
\begin{enumerate}
 \item The boundary of $C^+$ meets the boundary of $B^+$.
 \item Either the boundary of $C^+$ is a subset of $B^+$ or vice-versa.
\end{enumerate}

In the first case, let $v'$ be a point common to the two boundaries.
Then $v'$ lies on the boundary of $B^+$ and $v''$ lies on the boundary of $B^-$.
But these are subsets of the boundaries of $F^{\pm}$, and so $v$ is the
convex combination of at most 4 vertices.

In the second case, suppose that the boundary of $B^+$ is a subset of $C^+$.
Since $B^+$ is the intersection $F^+$ with a hyperplane, its boundary
must contain a vertex of $F^+$, as the set of vertices of $F^+$ is a connected
1-dimensional set whose convex hull is $F^+$.
Suppose that $v'\in B^+$ is a vertex of $F^+$.
Since $v''\in F^-$ is a convex combination of three vertices of
$F^-$, we see that $v$ is a convex combination of 1+3=4 vertices.
If the boundary of $C^+$ is a subset of $B^+$, then we may choose the point $v'$ in the
boundary of $C^+$ so that the corresponding point $v''$ is a a vertex of $F^-$.
Again, $v$ is the convex combination of $3+1=4$ vertices.
\QED


We wish to thank Ivano Bertini, Claudio Luchinat
and Giacomo Parigi of the Center for Magnetic Resonance of the
University of Florence for suggesting and discussing with us this
interesting problem.


\providecommand{\bysame}{\leavevmode\hbox to3em{\hrulefill}\thinspace}
\providecommand{\MR}{\relax\ifhmode\unskip\space\fi MR }
\providecommand{\MRhref}[2]{%
  \href{http://www.ams.org/mathscinet-getitem?mr=#1}{#2}
}
\providecommand{\href}[2]{#2}


\end{document}